\documentclass[11pt,a4paper]{article}
\usepackage{amssymb,amsmath}
\newcommand{\trp}{{\sf\scriptsize T}}
\newcommand\eop{\hfill$\Box$\\}
\newcommand\nst{n_{\scriptsize\rm st}}
\newcommand\nsts{n_{\tiny\rm st}}
\newtheorem{theorem}{Theorem}[section]
\newtheorem{lemma}[theorem]{Lemma}

\newcommand\asto{\,\stackrel{\rm\scriptsize a.s.}{\longrightarrow}\,}
\newcommand\dto{\,\stackrel{\rm\scriptsize d}{\longrightarrow}\,}
\newcommand\Ifkt {\hskip.1em 1\hskip-.6em 1\hskip.1em}
\newcommand\thb{\overline{\theta}}
\newcommand\LSn{\widehat{\theta}_n^{\rm\scriptsize (LS)}}

\begin{document}
\title{Convergence of least squares estimators
in the adaptive Wynn algorithm for a class of nonlinear regression models}
\author{Fritjof Freise$^1$, Norbert Gaffke$^2$, \ and \  Rainer Schwabe$^2$ \\[2ex]
\phantom{xxxxxxxxxx}\parbox{11cm}{\small $^1$TU Dortmund University \ and \ 
$^2$University of Magdeburg}
}
\date{}
\maketitle

\begin{abstract}
The paper continues the authors' work \cite{FF-NG-RS-18} on the adaptive Wynn algorithm
in a nonlinear regression model.  
In the present paper it is shown that if the mean response function satisfies a condition of  
`saturated identifiability', which was introduced by Pronzato \cite{Pronzato}, then the adaptive least squares estimators 
are strongly consistent. The condition states that
the regression parameter is identifiable under any saturated design, i.e.,
the values of the mean response function at any $p$ distinct design points determine the parameter point uniquely
where, typically, $p$ is the dimension of the regression parameter vector.    
Further essential assumptions are compactness of the experimental region and of the parameter space 
together with some natural continuity assumptions. 
If the true parameter point  is an interior point of the parameter space then
under some smoothness assumptions and asymptotic homoscedasticity of random errors  
the asymptotic normality of adaptive least squares estimators is obtained.
\end{abstract} 

\section{Introduction}
\setcounter{equation}{0}
The classical algorithm of Wynn \cite{Wynn} for D-optimal design in linear regression models  
has motivated a particular scheme for sequential adaptive design in nonlinear regression models,
see Freise \cite{Fritjof}, Pronzato \cite{Pronzato},  Freise, Gaffke, and Schwabe \cite{FF-NG-RS-18}.
We refer to this scheme as `the adaptive Wynn algorithm'. 
In our previous paper \cite{FF-NG-RS-18} the important class of 
generalized linear models with univariate response was considered, and the asymptotics of the sequences of designs
and maximum likelihood estimators under the adaptive Wynn algorithm was studied. In the present paper we
focus on another class of nonlinear models satisfying a condition introduced by Pronzato \cite{Pronzato}. We refer to that condition 
as `saturated identifiability'.  Under the adaptive Wynn algorithm the asymptotics of least squares estimators is studied.
As a main result, their strong consistency is obtained and, as a consequence, the asymptotic  
behavior of the generated design sequence becomes evident. Moreover, asymptotic normality of the adaptive least squares estimators
is obtained.   
      
Suppose a nonlinear regression model with mean response $\mu(x,\theta)$, $x\in{\cal X}$, $\theta\in\Theta$,
where ${\cal X}$ and $\Theta$ are the experimental region and the parameter space, resp.
Suppose that a family of $\mathbb{R}^p$-valued functions $f_\theta$, $\theta\in\Theta$, defined on ${\cal X}$
has been identified such that the $p\times p$ matrix $f_\theta(x)\,f_\theta^\trp(x)$ is the elementary 
information matrix of $x\in{\cal X}$ at $\theta\in\Theta$. Note that a vector $a\in\mathbb{R}^p$ is
written as a column vector and $a^\trp$ denotes its transposed which is thus a $p$-dimensional row vector.
A design $\xi$ is a probability measure on ${\cal X}$ with finite support. That is, $\xi$
is described by its support, denoted by ${\rm supp}(\xi)$, which is a nonempty finite subset of ${\cal X}$,
and by its weights $\xi(x)$ for $x\in{\rm supp}(\xi)$ which are positive real numbers with 
$\sum_{x\in {\scriptsize\rm supp}(\xi)}\xi(x)\,=1$.
The information matrix of a design $\xi$ at $\theta\in\Theta$ is defined by
\begin{equation}
M(\xi,\theta)\,=\,\sum_{x\in {\scriptsize\rm supp}(\xi)}\xi(x)\,f_\theta(x)\,f_\theta^\trp(x),
\label{eq1-1}
\end{equation} 
which is a nonnegative definite $p\times p$ matrix. 

In applications the family $f_\theta$, $\theta\in\Theta$, will be related to the mean response $\mu(x,\theta)$,
$x\in{\cal X}$, $\theta\in\Theta$.
In particular, $f_\theta(x)$ may be the gradient of $\mu(x,\theta)$ w.r.t.~$\theta$ for each fixed $x$ where
$\Theta\subseteq\mathbb{R}^p$. However, we do not generally assume a relation 
between $\mu$ and the family $f_\theta$, $\theta\in\Theta$. 
An exception is a result on the asymptotic normality of the least squares estimators
in Section 3 in which case the assumed relation is explicitly stated. 
Throughout we assume the following basic conditions (B1)-(B4).
\begin{itemize}
\item[(B1)]
The experimental region ${\cal X}$ is a compact metric space.
\item[(B2)] 
The parameter space $\Theta$ is a compact metric space.
\item[(B3)] 
The real-valued mean response function \ $(x,\theta)\mapsto \mu(x,\theta)$, defined on the Cartesian product space  
${\cal X}\times \Theta$, is continuous.
\item[(B4)]
The family $f_\theta$, $\theta\in\Theta$, of $\mathbb{R}^p$-valued functions on ${\cal X}$ satisfies:\\
(i) for each $\theta\in\Theta$ the image  $f_\theta({\cal X})$ \ spans $\mathbb{R}^p$;\\
(ii) the function \ $(x,\theta)\mapsto f_\theta(x)$, defined on the Cartesian product space  
${\cal X}\times \Theta$, is continuous.
\end{itemize}
By $\mathbb{N}$ and $\mathbb{N}_0$ we denote the set of all positive integers and all nonnegative integers, resp.
By $\delta_x$ for any $x\in{\cal X}$ 
we denote the one-point probability distribution on ${\cal X}$ concentrated at the point $x$.
The adaptive Wynn algorithm collects iteratively design points $x_i\in{\cal X}$, $i\in\mathbb{N}$, while 
adaptively estimating $\theta$ on the basis of the current design points and observed responses at those points.  
In greater detail the algorithm reads as follows. \\[.5ex]
\underbar{Adaptive Wynn algorithm.}
\begin{itemize}
\item[(o)] {\em Initialization:}    
A suitable number $\nst\in\mathbb{N}$ and design points $x_1,\ldots,x_{\nsts}\in{\cal X}$ are chosen such that
the starting design $\xi_{\nsts}=\frac{1}{\nst}\sum_{i=1}^{\nsts} \delta_{x_i}$ has positive definite information matrices, 
i.e., for all $\theta\in\Theta$ the information matrix $M(\xi_{\nsts},\theta)$ is positive definite.
Observed responses $y_1,\ldots,y_{\nsts}$ at the design points $x_1,\ldots,x_{\nsts}$ are taken, and 
based on those an initial parameter estimate $\theta_{\nsts}$ is calculated,
\[
\theta_{\nsts}=\widehat{\theta}_{\nsts}(x_1,y_1,\ldots,x_{\nsts},y_{\nsts})\,\in\Theta.
\]
\item[(i)] {\em Iteration:} 
At stage $n\ge\nst$ the current data is given by the points $x_1,\ldots,x_n\in{\cal X}$ which form the design
$\xi_n=\frac{1}{n}\sum_{i=1}^n\delta_{x_i}$, and by the observed responses  $y_1,\ldots,y_n$ at $x_1,\ldots,x_n$, resp.,
along with a parameter estimate $\theta_n$,
\begin{equation}
\theta_n=\widehat{\theta}_n(x_1,y_1,\ldots,x_n,y_n)\,\in\Theta.\label{eq1-0}
\end{equation}
The iteration rule is given by
\begin{equation}
x_{n+1}\,=\,\arg\ \max_{x\in{\cal X}}f_{\theta_n}^\trp(x)\,M^{-1}(\xi_n,\theta_n)\,
 f_{\theta_n}(x). \label{eq1-1}
\end{equation}  
An observation $y_{n+1}$ of the response at $x_{n+1}$ is taken and a new parameter estimate 
$\theta_{n+1}$ based on the augmented data is computed,
\[
\theta_{n+1}=\widehat{\theta}_{n+1}(x_1,y_1,\ldots,x_n,y_n,x_{n+1},y_{n+1})\,\in\Theta.
\]
Replace $n$ by $n+1$ and repeat the iteration step (i).
\end{itemize}
\eop
Of course, (\ref{eq1-1}) requires the information matrix $M(\xi_n,\theta_n)$ to be positive definite at each stage $n\ge\nst$.
In fact, this is ensured by the choice of the initial design $\xi_{\nsts}$ since, obviously,  
the sequence of designs $\xi_n$, $n\ge\nst$ satisfies
\begin{eqnarray}
&&\textstyle\xi_{n+1}\,=\,\frac{n}{n+1}\,\xi_n\,+\,\frac{1}{n+1}\,\delta_{x_{n+1}},\ \ n\ge\nst, \label{eq1-2}\\
&&\textstyle M(\xi_{n+1},\theta)\,=\,\frac{n}{n+1}\,M(\xi_n,\theta)\,+\,\frac{1}{n+1}\,
f_\theta(x_{n+1})\,f_\theta^\trp(x_{n+1}), \ \ n\ge\nst,\  \theta\in\Theta,\label{eq1-3}
\end{eqnarray}
from which one concludes by induction that 
$M(\xi_n,\theta)$ is positive definite for all $n\ge\nst$ and all $\theta\in\Theta$. 
The existence of an initial design $\xi_{\nsts}$ as required will be shown in Section 2.
 
The algorithm uses, in particular, an observed response $y_i$ at each current design point $x_i$. 
So the generated sequence of design points, $x_i$, $i\in\mathbb{N}$, and the
corresponding sequence of designs $\xi_n$, $n\ge\nst$, are random sequences with a particular dependence structure
caused by (\ref{eq1-0}) and (\ref{eq1-1}). An appropriate stochastic model  
will be stated in Section 3 which was used in
Freise, Gaffke, and Schwabe \cite{FF-NG-RS-18} and goes back to Lai and Wei \cite{Lai-Wei}, Lai \cite{Lai}, and     
Chen, Hu, and Ying \cite{Chen-Hu-Ying}. In particular, the generated sequence $x_i$, $i\in\mathbb{N}$,
and the observed responses $y_i$, are viewed as values of random variables $X_i$, $i\in\mathbb{N}$, and
$Y_i$, $i\in\mathbb{N}$, resp., following a stochastic model which we call an `adaptive regression model'.
Our formulation of the adaptive Wynn algorithm is a description of the paths of the 
stochastic sequence $(X_i,Y_i)$, $i\in\mathbb{N}$.     

The estimators $\widehat{\theta}_n$, $n\ge\nst$, employed by the algorithm to produce the
estimates $\theta_n$, $n\ge\nst$, in (\ref{eq1-0}), in principle, may be any estimators of $\theta$ such that their values are 
in $\Theta$ and $\widehat{\theta}_n$ is a function of the data $x_1,y_1,\ldots,x_n,y_n$ available at stage $n$.
Such estimators will be called adaptive estimators.
Later, strong consistency of $\widehat{\theta}_n$, $n\ge\nst$, will be required. In Section 3 we focus on adaptive
least squares estimators $\widehat{\theta}_n^{({\rm\scriptsize LS})}$, i.e.,  
\[
\widehat{\theta}_n^{({\rm\scriptsize LS})}(x_1,y_1,\ldots,x_n,y_n)\,=\,\arg\min_{\theta\in\Theta}  
\sum_{i=1}^n\bigl(y_i-\mu(x_i,\theta)\bigr)^2.
\]
Note that when dealing with the adaptive least squares estimators we will not necessarily assume that 
the estimators $\widehat{\theta}_n$, $n\ge\nst$, employed by the algorithm 
are given by the adaptive least squares estimators. For proving strong consistency of the latter, irrespective
which adaptive estimators $\widehat{\theta}_n$ are used in the algorithm, we need 
an additional condition on the mean response $\mu$ which we call `saturated identifiability' (SI),
and which was introduced by Pronzato \cite{Pronzato} in the case of a finite experimental region ${\cal X}$.
\begin{itemize}
\item[(SI)]
If $z_1,\ldots,z_p\in{\cal X}$ are pairwise distinct design points then the $\mathbb{R}^p$-valued function
on $\Theta$,
\[
\theta\longmapsto\bigl(\mu(z_1,\theta),\ldots,\mu(z_p,\theta)\bigr)^\trp,
\]
is an injection, i.e., if $\theta,\theta'\in\Theta$ and $\mu(z_j,\theta)=\mu(z_j,\theta')$, $1\le j\le p$, 
then $\theta=\theta'$.
\end{itemize}
The employed stochastic model for the adaptive Wynn algorithm 
includes a martingale difference scheme for the error variables. So limit theorems for martingales  
can be applied: A Strong Law of Large Numbers and a Central Limit Theorem to prove strong consistency and
asymptotic normality, resp., of the adaptive least squares estimators (Theorem 3.1 and Theorem 3.2).   
As a remark on Pronzato \cite{Pronzato} we note that, in our view, his proof of Theorem 1 of that paper (pp.~210-211) 
is not convincing when applying the Law of Iterated Logarithm to {\em random} subsequences of the error variables.  
Arguments by martingales might help, like those in Pronzato \cite{Pronzato-b}.     

We start with some auxiliary results which are the content of Section 2. 

\section{Auxiliary results}
\setcounter{equation}{0}
Throughout we assume (B1)-(B4) as introduced in the previous section. 
Note, however, that (B3) will not play a role in this section. 
Firstly, we give a proof of the existence of an initial design as required      
in the algorithm.

\begin{lemma}\quad\label{lem2-1}\\
There exist an $\nst\in\mathbb{N}$ and design points 
$x_1,\ldots,x_{\nsts}\in{\cal X}$ such that for every $\theta\in\Theta$ the vectors
\ $f_\theta(x_1),\ldots,f_\theta(x_{\nsts})$ \ span $\mathbb{R}^p$. Hence, for such $x_i$, $1\le i\le \nst$,
the design \ $\xi_{\nsts}=\frac{1}{\nst}\sum_{i=1}^{\nsts}\delta_{x_i}$ has the property that its information matrix
$M(\xi_{\nsts},\theta)$ is positive definite for all $\theta\in\Theta$. 
\end{lemma}

\noindent
{\bf Proof.} \ By (B4)-(i), for each $\theta\in\Theta$ there exist $p$ design points 
$z_1(\theta),\ldots,z_p(\theta)\in{\cal X}$ such that the vectors $f_\theta\bigl(z_1(\theta)\bigr),\ldots,
f_\theta\bigl(z_p(\theta)\bigr)$ are linearly independent. By (B2) and (B4) (ii), for each $\theta\in\Theta$
the set
\[
U(\theta)\,=\,\Bigl\{\tau\in\Theta\,:\,\det\,\Bigl[f_\tau\bigl(z_1(\theta)\bigr),\ldots,
f_\tau\bigl(z_p(\theta)\bigr)\Bigr]\not=0\Bigr\}
\]
is an open set in the (compact) metric space $\Theta$, and
$\theta\in U(\theta)$. Hence, trivially, 
$\Theta=\bigcup_{\theta\in\Theta}U(\theta)$, and by (B2) there is an $r\in\mathbb{N}$
and points $\theta_1,\ldots,\theta_r\in\Theta$ such that \ $\Theta=\bigcup_{j=1}^rU(\theta_j)$.
Denote $x_{ij}=z_i(\theta_j)$, $1\le i\le p$, $1\le j\le r$. Then, for every $\tau\in\Theta$ 
the set of vectors 
\[
\bigl\{f_\tau(x_{ij})\,:\,1\le i\le p,\ 1\le j\le r\bigr\}
\]
spans $\mathbb{R}^p$. In fact, for any given $\tau\in\Theta$ there is some $j_0\in\{1,\ldots,r\}$ with $\tau\in U(\theta_{j_0})$ hence
$\det\,\bigl[f_\tau(x_{1j_0}),\ldots,f_\tau(x_{pj_0})\bigr]\not=0$, i.e.,
the vectors $f_\tau(x_{1j_0}),\ldots,f_\tau(x_{pj_0})$ constitute a basis of $\mathbb{R}^p$. 
So, for $\nst=pr$ and $x_1,\ldots,x_{\nsts}$ being a relabelled family of the points $x_{ij}$, $1\le i\le p$, $1\le j\le r$,
the vectors $f_\tau(x_1),\ldots,f_\tau(x_{\nsts})$ span $\mathbb{R}^p$, and hence the information matrix
$M(\xi_{\nsts},\tau)=\frac{1}{\nst}\sum_{i=1}^{\nsts}f_\tau(x_i)\,f^\trp_\tau(x_i)$ is positive definite.
\eop

Let any path of the adaptive Wynn algorithm be given as described in the previous section.
In particular, $x_i$, $i\in\mathbb{N}$, is the sequence of design points and 
$\xi_n=\frac{1}{n}\sum_{i=1}^n\delta_{x_i}$, $n\ge\nst$, is the corresponding sequence of designs.
For the following two lemmas 
no assumption on the employed adaptive estimators $\widehat{\theta}_n$, $n\ge\nst$, is needed.
In other words, the sequence $\theta_n$, $n\ge\nst$, of parameter estimates appearing in the path
may be arbitrary. 
Note that Lemma \ref{lem2-3} below extends Lemma 2 of Pronzato \cite{Pronzato} 
who restricted to a finite experimental region.

We denote the distance function in the compact metric space ${\cal X}$  by ${\rm d}_{\cal X}(x,z)$, $x,z\in{\cal X}$.
If $S_1$ and $S_2$ are nonempty subsets of ${\cal X}$ then the distance ${\rm d}_{\cal X}(S_1,S_2)$ of $S_1$ and $S_2$ is defined by
${\rm d}_{\cal X}(S_1,S_2)\,=\,\inf\{{\rm d}_{\cal X}(x,z)\,:\,x\in S_1,\ z\in S_2\}$. In case that $S_1=\{x\}$ is a singleton we write
${\rm d}_{\cal X}(x,S_2)$ instead of ${\rm d}_{\cal X}(\{x\},S_2)$.   
If $S$ is a nonempty subset of ${\cal X}$ then the diameter of $S$
is defined by ${\rm diam}(S)\,=\,\sup\{{\rm d}_{\cal X}(x,z)\,:\,x,z\in S\}$. 

\begin{lemma}\quad\label{lem2-2}\\
Suppose $p\ge2$. Let $\varepsilon>0$ be given. 
Then there exist $d>0$ and  $n_0\ge\nst$ such that 
\[
\xi_n(S)\le\frac{1}{p}+\varepsilon\ 
\mbox{ for all $\emptyset\not=S\subseteq{\cal X}$ with ${\rm diam}(S)\le d$ and all $n\ge n_0$.}
\]
\end{lemma}  
       
\noindent
{\bf Proof.} \ 
Without loss of generality we may assume  $\varepsilon<1-p^{-1}$.
In \cite{FF-NG-RS-18} we introduced the positive real constants
\begin{equation}
\gamma =\sup_{x\in{\cal X},\, \theta\in\Theta}\Vert f_\theta(x)\Vert\ \quad\mbox{and}\quad
\kappa = \inf_{\Vert v\Vert=1,\, \theta\in\Theta}\ \max_{x\in{\cal X}}\bigl(v^\trp f_\theta(x)\bigr)^2,
\label{eq2-1}
\end{equation}
and Lemma 2.3 of that paper stated the following. 
\begin{eqnarray}
&&\mbox{If $0<\eta<1-p^{-1/2}$, $n\ge \nst$, and $S\subseteq{\cal X}$
are given such that}\nonumber\\
&&\mbox{$\Vert f_{\theta_n}(x)-f_{\theta_n}(z)\Vert \le\eta\kappa/\gamma$ for all $x,z\in S$
and $\xi_n(S)>(1-\eta)^{-2}p^{-1}$,}\nonumber\\ 
&&\mbox{then \ $x_{n+1}\not\in S$.}\label{eq2-2}
\end{eqnarray} 
Choose \ $\eta:=1-\bigl(1+p\varepsilon/2\bigr)^{-1/2}$. Then 
$0<\eta<1-p^{-1/2}$ and $(1-\eta)^{-2}p^{-1}=p^{-1}+\varepsilon/2$. 
By (B1), (B2), and (B4) the function \ $(x,\theta)\mapsto f_\theta(x)$ \ is uniformly continuous on its compact   
domain ${\cal X}\times \Theta$. So there exists a $d>0$ such that
\begin{equation}
\mbox{if $x,z\in{\cal X}$ and ${\rm d}_{\cal X}(x,z)\le d$ then \ $\Vert f_\theta(x)-f_\theta(z)\Vert\le\eta\kappa/\gamma$
\ $\forall\  \theta\in\Theta$.} \label{eq2-3}
\end{equation}
We show that $d$ fulfills the requirement of the lemma. 
Let $\emptyset\not=S\subseteq{\cal X}$ with ${\rm diam}(S)\le d$. By (\ref{eq2-3}) and (\ref{eq2-2})
the sequence $\xi_n(S)$, $n\ge \nst$, has the property that for all $n\ge \nst$,
\begin{eqnarray*}
\xi_{n+1}(S)& =&  \textstyle\frac{n}{n+1}\xi_n(S)\quad \mbox{ if }\xi_n(S)>\frac{1}{p}+\frac{\varepsilon}{2},\\
\xi_{n+1}(S)& \le &\textstyle \xi_n(S)+\frac{1}{n+1} \  \mbox{ if }\xi_n(S)\le\frac{1}{p}+\frac{\varepsilon}{2}. 
\end{eqnarray*}
An application of Lemma 2.1 in \cite{FF-NG-RS-18} to the sequence $\beta_n:=\xi_n(S)$, $n\ge\nst$, and
$\beta:=\frac{1}{p}+\frac{\varepsilon}{2}$, $\widetilde{\beta}:=\frac{1}{p}+\varepsilon$ yields that
\[
\xi_n(S)\le\frac{1}{p}+\varepsilon \ \mbox{ for all }n\ge 
n_0:=\big\lceil{\textstyle(\frac{1}{p}+\frac{\varepsilon}{2})^{-1}}\big\rceil\cdot
\max\bigl\{\nst,\lceil 2/\varepsilon\rceil\bigr\},
\]
where $\lceil a\rceil$, for $a\in\mathbb{R}$, denotes the smallest integer greater than or equal to $a$. 
Since $n_0$ does not depend on the particular set $S$ the result follows.
\eop
 
\begin{lemma}\quad\label{lem2-3}\\
Suppose $p\ge2$. There exist $n_0\ge\nst$, $\pi_0>0$, and $d_0>0$ such that the following holds.\\[1ex]
\indent For each $n\ge n_0$ there are $p$ subsets $S_{1,n},S_{2,n},\ldots,S_{p,n}$ of ${\cal X}$ such that\\
$\xi_n(S_{j,n})\ge \pi_0$, $1\le j\le p$, \ 
${\rm diam}(S_{j,n})\le d_0$, $1\le j\le p$, and \ ${\rm d}_{\cal X}(S_{j,n},S_{k,n})\ge d_0$,  
$1\le j< k\le p$.
\end{lemma} 

\noindent
{\bf Proof.} \ 
Fix an $\varepsilon$ with $0<\varepsilon<\frac{1}{p(p-1)}$. 
Choose  $d>0$ and  $n_0\ge \nst$ according to Lemma \ref{lem2-2}. 
By compactness of ${\cal X}$ there is a positive integer $q$
and nonempty subsets $R_1,\ldots,R_q$ of ${\cal X}$ such that 
\[
{\cal X}\,=\,\bigcup_{\ell=1}^qR_\ell\ \mbox{ and }\ {\rm diam}(R_\ell)\le d/3\ \mbox{for all $\ell=1,\ldots,q$.}
\]
We show that \  $n_0$, $\pi_0:=\bigl(p^{-1}-(p-1)\varepsilon\bigr)/q$, and $d_0:=d/3$ \ satisfy the requirements of the 
assertion. To this end let $n\ge n_0$ be given. We construct inductively subsets $S_{j,n}$, $1\le j\le p$, as required.\\
\underbar{$j=1$}\,: Clearly, $\sum_{\ell=1}^q\xi_n(R_\ell)\ge1$. Choose $\ell_n\in\{1,\ldots,q\}$ achieving the maximum value of $\xi_n(R_\ell)$,
$1\le \ell\le q$,
and set \ $S_{1,n}:=R_{\ell_n}$. Then 
$\xi_n(S_{1,n})=\max_{1\le \ell\le q}\xi_n(R_\ell)\ge1/q\ge \pi_0$ and ${\rm diam}(S_{1,n})={\rm diam}(R_{\ell_n})\le d_0$.\\
\underbar{Induction step}: Let an $r\in\{1,\ldots,p-1\}$ be given along with subsets $S_{1,n},\ldots,S_{r,n}$ of ${\cal X}$ such that
$\xi_n(S_{j,n})\ge \pi_0$ and ${\rm diam}(S_{j,n})\le d_0$,  $1\le j\le r$, and ${\rm d}_{\cal X}(S_{j,n},S_{k,n})\ge d_0$, $1\le j<k\le r$. 
Let $\overline{S}_{j,n}:=\bigl\{x\in{\cal X}\,:\,{\rm d}_{\cal X}(x,S_{j,n})\le d_0\bigr\}$, $1\le j\le r$. As it is easily seen,
${\rm diam}(\overline{S}_{j,n})\le3d_0=d$ and hence $\xi_n(\overline{S}_{j,n})\le\frac{1}{p}+\varepsilon$. 
So for \ $T_{r,n}:=\bigcup_{j=1}^r\overline{S}_{j,n}$ \ one has \ $\xi_n(T_{r,n})\le r\,\bigl(\frac{1}{p}+\varepsilon\bigr)$,
\ and hence  
\[
\xi_n({\cal X}\setminus T_{r,n})\ge 1- r\,\bigl(\frac{1}{p}+\varepsilon\bigr)\ge 1- (p-1)(\frac{1}{p}+\varepsilon)=p^{-1}-(p-1)\varepsilon.
\]
Observing that \ ${\cal X}\setminus T_{r,n}=\bigcup_{\ell=1}^q(R_\ell\setminus T_{r,n})$ one gets
\[
p^{-1}-(p-1)\varepsilon\,\le \sum_{\ell=1}^q \xi_n(R_\ell\setminus T_{r,n}).
\]
Choose $\ell_n\in\{1,\ldots,q\}$ which achieves the maximum value of $\xi_n(R_\ell\setminus T_{r,n})$, $1\le \ell\le q$, 
and set $S_{r+1,n}:= R_{\ell_n}\setminus T_{r,n}$. Then  
$\xi_n(S_{r+1,n})=\max_{1\le \ell\le q}\xi_n\bigl(R_\ell\setminus T_{r,n}\bigr)\ge \bigl(p^{-1}-(p-1)\varepsilon\bigr)/q=\pi_0$ and
${\rm diam}(S_{r+1,n})\le{\rm diam}(R_{\ell_n})\le d_0$. Moreover for each $j=1,\ldots,r$, since $S_{r+1,n}\cap \overline{S}_{j,n}=\emptyset$,
one has ${\rm d}_{\cal X}(x,S_{j,n})>d_0$ for all $x\in S_{r+1,n}$ and hence  ${\rm d}_{\cal X}(S_{r+1,n},S_{j,n})\ge d_0$. So 
we have subsets $S_{1,n},\ldots,S_{r,n},S_{r+1,n}$ such that
\begin{eqnarray*}
&&\xi_n(S_{j,n})\ge \pi_0\ \mbox{ and }\ {\rm diam}(S_{j,n})\le d_0,\ \ 1\le j\le r+1,\\
&&{\rm d}_{\cal X}(S_{j,n},S_{k,n})\ge d_0,\ \ 1\le j<k\le r+1. 
\end{eqnarray*}
This completes the inductive construction and the proof of the lemma.
\eop

\vspace*{1ex}\noindent
\underbar{Remark.} \ In the case that ${\cal X}$ is finite it is easily seen that 
in Lemma \ref{lem2-3} the subsets $S_{1,n},\ldots,S_{p,n}$ can be chosen to be singletons
for all $n\ge n_0$.
So, in this case, the lemma yields the result of Lemma 2 of Pronzato \cite{Pronzato}.
\eop

\section{Convergence of least squares estimators} 
\setcounter{equation}{0}
For an analysis of the adaptive Wynn algorithm, the generated sequence $x_i$, $i\in\mathbb{N}$,
and the observed responses $y_i$, are  viewed as values of random variables $X_i$, $i\in\mathbb{N}$, and
$Y_i$, $i\in\mathbb{N}$, resp., whose dependence structure is described by the following two
assumptions (A1) and (A2), see \cite{FF-NG-RS-18}, \cite{Lai-Wei}, \cite{Lai}, \cite{Chen-Hu-Ying}. 
The model thereby stated might be called an `adaptive regression model'.
By $\overline{\theta}$ we denote the true point of the parameter space $\Theta$ governing the data. 
All the random variables appearing in this section are thought to be defined on a common probability   
space $(\Omega,{\cal F},\mathbb{P}_{\overline{\theta}})$, where $\Omega$ is a nonempty set, ${\cal F}$
is a sigma-field of subsets of $\Omega$, and $\mathbb{P}_{\overline{\theta}}$ is a probability measure on ${\cal F}$ 
corresponding to the true parameter point $\overline{\theta}$. We assume, as before, the conditions (B1)-(B4),
and now additionally the following conditions (A1) and (A2) constituting the adaptive regression model.
\begin{itemize}
\item[(A1)]
There is given a nondecreasing sequence of sub-sigma-fields of ${\cal F}$, 
\ ${\cal F}_0\subseteq{\cal F}_1\subseteq\,\ldots\,\subseteq{\cal F}_n\subseteq\,\ldots$ 
such that  for each $i\in\mathbb{N}$ the random variable $X_i$ is ${\cal F}_{i-1}$-measurable
and the random variable $Y_i$ is ${\cal F}_i$-measurable.
\item[(A2)]
$Y_i\,=\,\mu(X_i,\overline{\theta})\,+\,e_i$ \ with real-valued 
square integrable random errors $e_i$ such that\\ 
${\rm E}\bigl(e_i\,\big|\,{\cal F}_{i-1}\bigr)\,=0 \ \mbox{\,a.s.}$\  for all $i\in\mathbb{N}$, and \ 
$\sup_{i\in\mathbb{N}}{\rm E}\bigl(e_i^2\,\big|\,{\cal F}_{i-1}\bigr)\,<\infty \ \ \mbox{a.s.}$
\end{itemize}
As before, $\widehat{\theta}_n$, $n\ge\nst$, are the adaptive estimators employed by the algorithm,
now viewed as random variables, $\widehat{\theta}_n=\widehat{\theta}_n(X_1,Y_1,\ldots,X_n,Y_n)$.
Of course, a desirable property of these estimators would be strong consistency,  
i.e., almost sure convergence to $\overline{\theta}$ (as $n\to\infty$), for short
$\widehat{\theta}_n\asto\overline{\theta}$.
As shown in our previous paper \cite{FF-NG-RS-18}, Corollary 3.2, if the estimators $\widehat{\theta}_n$
are strongly consistent,  then the sequence $\xi_n$, $n\ge\nst$, of (random) designs generated by the algorithm
is almost surely asymptotically D-optimal in the sense that \ 
$M(\xi_n,\widehat{\theta}_n)\asto M(\xi^*_{\overline{\theta}},\overline{\theta})$, where   
$\xi^*_{\overline{\theta}}$ is a locally D-optimal design at $\overline{\theta}$. 
In fact, the conclusion of that corollary is stronger: if  the estimators $\widehat{\theta}_n$
are strongly consistent then $M(\xi_n,\widetilde{\theta}_n)\asto M(\xi^*_{\overline{\theta}},\overline{\theta})$
holds for {\em every} strongly consistent sequence of $\Theta$-valued estimators $\widetilde{\theta}_n$.  
Under the condition (SI) of saturated identifyability introduced in Section 1,
the next result shows that strong consistency of $\widehat{\theta}_n$ holds when employing the least squares estimators,  
i.e., $\widehat{\theta}_n=\widehat{\theta}_n^{(\rm\scriptsize LS)}$, $n\ge\nst$.
Actually, the result is stronger: when employing any adaptive estimators  
$\widehat{\theta}_n$, $n\ge\nst$, the sequence of least squares estimators  
$\widehat{\theta}_n^{(\rm\scriptsize LS)}$, $n\ge\nst$, is strongly consistent (under condition (SI)). 

\begin{theorem}\quad\label{theo3-1}\\
Assume condition (SI). Then, irrespective of the employed sequence of adaptive estimators $\widehat{\theta}_n$ 
in the algorithm, the sequence of adaptive least squares estimators 
$\widehat{\theta}_n^{(\rm\scriptsize LS)}$  is strongly consistent: \ 
$\widehat{\theta}_n^{(\rm\scriptsize LS)}\asto\overline{\theta}$.
\end{theorem}

\vspace*{1ex}\noindent
{\bf Proof.} \ Define for all $n\in\mathbb{N}$ and $\theta\in\Theta$  random variables
\[
S_n(\theta):=\sum_{i=1}^n\bigl(Y_i-\mu(X_i,\theta)\bigr)^2\quad \mbox{and}\quad 
D_n(\theta,\overline{\theta})\,:=\,\sum_{i=1}^n\bigl(\mu(X_i,\theta)-\mu(X_i,\overline{\theta})\bigr)^2.
\]
The proof is divided into three steps. For $\varepsilon>0$ we denote 
$C(\overline{\theta},\varepsilon):=
\bigl\{\theta\in\Theta\,:\,{\rm d}_\Theta(\theta,\overline{\theta})\ge\varepsilon\bigr\}$,
where ${\rm d}_\Theta$ denotes the distance function in $\Theta$.\\
\underbar{Step 1.} Show that for all $\varepsilon>0$ with 
$C(\overline{\theta},\varepsilon)\not=\emptyset$,
\[
\Big|\,\frac{1}{n}\Bigl(\inf_{\theta\in C(\overline{\theta},\varepsilon)} S_n(\theta)\,-S_n(\overline{\theta})\Bigr)\,
-\,\frac{1}{n}\inf_{\theta\in C(\overline{\theta},\varepsilon)} D_n(\theta,\overline{\theta})\,\Big|
\ \asto\,0.
\]
\underbar{Step 2.} Show that for all $\varepsilon>0$ with 
$C(\overline{\theta},\varepsilon)\not=\emptyset$,
\[
\liminf_{n\to\infty}\Bigl(\frac{1}{n}\,
\inf_{\theta\in C(\overline{\theta},\varepsilon)}D_n(\theta,\overline{\theta})\Bigr)\ >\,0\ \ \mbox{a.s.}
\]
\underbar{Step 3.} Conclude from the results of Step 1 and Step 2 that for all $\varepsilon>0$ with 
$C(\overline{\theta},\varepsilon)\not=\emptyset$, 
\begin{equation}
\inf_{\theta\in C(\overline{\theta},\varepsilon)}S_n(\theta)\,-\,S_n(\overline{\theta})\ \asto\,\infty.
\label{eq3-1}
\end{equation}
From (\ref{eq3-1}), applying Lemma 1 of Wu \cite{Wu}, one gets \ 
$\widehat{\theta}_n^{(\rm\scriptsize LS)}\asto\overline{\theta}$.\\
\underbar{Ad Step 1.} As in Pronzato \cite{Pronzato}, p.~230, one calculates
\begin{eqnarray*}
&&S_n(\theta)-S_n(\overline{\theta})\,=\,D_n(\theta,\overline{\theta})\,+\,2W_n(\theta,\overline{\theta}),\ \ \mbox{where}\\
&&W_n(\theta,\overline{\theta})\,:=\,\sum_{i=1}^n
\bigl(\mu(X_i,\overline{\theta})-\mu(X_i,\theta)\bigr)\,e_i.
\end{eqnarray*}
It follows that
\[
\Big|\,\frac{1}{n}\Bigl(\inf_{\theta\in C(\overline{\theta},\varepsilon)} S_n(\theta)\,-S_n(\overline{\theta})\Bigr)\,
-\,\frac{1}{n}\inf_{\theta\in C(\overline{\theta},\varepsilon)} D_n(\theta,\overline{\theta})\,\Big|
\le\ \frac{2}{n}\,\sup_{\theta\in\Theta}\big|W_n(\theta,\overline{\theta})\bigr|.
\]
Applying Lemma 3.1, part (c), in \cite{FF-NG-RS-18} with $h(x,\theta)=\mu(x,\overline{\theta})-\mu(x,\theta)$,
$(x,\theta)\in{\cal X}\times\Theta$, the result of Step 1 follows.\\
\underbar{Ad Step 2.}  Consider any path $x_i,y_i$, $i\in\mathbb{N}$, and $\theta_n$, $n\ge\nst$
of the  sequences $X_i,Y_i$, $i\in\mathbb{N}$,
and $\widehat{\theta}_n$, $n\ge\nst$. Firstly, consider the simple case $p=1$. Then condition (SI)
implies that $\mu(x,\theta)\not=\mu(x,\thb)$ for all $\theta\in C(\overline{\theta},\varepsilon)$,
and hence by (B3) 
\[
c_\varepsilon\,:=\,\inf_{x\in{\cal X}}\bigl(\mu(x,\theta)-\mu(x,\thb)\bigr)^2\ >0.
\]
It follows that $\frac{1}{n}\inf_{\theta\in C(\thb,\varepsilon)}D_n(\theta,\thb)\,\ge c_\varepsilon$ for all $n$
and, in particular,  its limit inferior is positive. Now let $p\ge2$. 
According to Lemma \ref{lem2-3}, choose $n_0\ge\nst$, $\pi_0>0$, $d_0>0$, and subsets
$S_{1,n},\ldots, S_{p,n}\subseteq{\cal X}$ for all $n\ge n_0$.  
Define a subset of the $p$-fold product space ${\cal X}^p$ by
\[
\Delta\,:=\,\bigl\{(z_1,\ldots,z_p)\in{\cal X}^p\,:\,{\rm d}_{\cal X}(z_j,z_k)\ge d_0,\ \ 1\le j<k\le p\bigr\}.
\]
By (SI), $\sum_{j=1}^p\bigl(\mu(z_j,\theta)-\mu(z_j,\overline{\theta})\bigr)^2>0$ for all $(z_1,\ldots,z_p)\in\Delta$
and all $\theta\not=\overline{\theta}$. By (B1) the set $\Delta$ is compact and by (B2) the set 
$C(\overline{\theta},\varepsilon)$ is compact. So, together with (B3), one concludes that the following infimum $c_\varepsilon$
is positive,
\[
c_\varepsilon:=\inf\Bigl\{\sum_{j=1}^p\bigl(\mu(z_j,\theta)-\mu(z_j,\overline{\theta})\bigr)^2\,:\,
   (z_1,\ldots,z_p)\in\Delta,\ \theta\in C(\overline{\theta},\varepsilon)\Bigr\}.
\]
For all $n\ge n_0$ and all permutations $\sigma$ of $\{1,\ldots,p\}$ the Cartesian product \ 
$S_n^\sigma:=S_{\sigma(1),n}\times S_{\sigma(2),n}\times\ldots\times S_{\sigma(p),n}$
is a subset of $\Delta$, hence \ $R_n:=\bigcup_\sigma S_n^\sigma\subseteq\Delta$.
Note that $S_n^\sigma\cap S_n^\tau=\emptyset$ for any two different permutations $\sigma$ and $\tau$.
Consider the $p$-fold product measure $\xi_n^p$. Then, for all $\sigma$ and all $n\ge n_0$ one has 
$\xi_n^p(S_n^\sigma)=\prod_{j=1}^p\xi_n(S_{\sigma(j),n})\ge \pi_0^p$ and hence $\xi_n(R_n)\ge p!\,\pi_0^p$.
So
\[
\int_{{\cal X}^p}\sum_{j=1}^p\bigl(\mu(z_j,\theta)-\mu(z_j,\overline{\theta})\bigr)^2\,{\rm d}\xi_n^p(z_1,\ldots,z_p)
\ \ge c_\varepsilon\,p!\,\pi_0^p \ \mbox{ for all $n\ge n_0$ and $\theta\in C(\overline{\theta},\varepsilon)$.}
\]
The integral on the l.h.s.\ of that inequality is equal to 
\[
p\,\int_{{\cal X}}\bigl(\mu(z,\theta)-\mu(z,\overline{\theta})\bigr)^2\,{\rm d}\xi_n(z)\,=\,\frac{p}{n}
\sum_{i=1}^n \bigl(\mu(x_i,\theta)-\mu(x_i,\overline{\theta})\bigr)^2.
\]
It follows that 
\[
\frac{1}{n}\,\inf_{\theta\in C(\overline{\theta},\varepsilon)}
\sum_{i=1}^n\bigl(\mu(x_i,\theta)-\mu(x_i,\overline{\theta})\bigr)^2\,\ge\,
c_\varepsilon\,(p-1)!\,\pi_0^p\ \ \forall\ n\ge n_0,
\]
which implies that the limit inferior of the l.h.s.\ of that inequality is positive.\\
\underbar{Ad Step 3.} By the results of Step 1 and Step 2,
\[
\liminf_{n\to\infty}\frac{1}{n}\Bigl(\inf_{\theta\in C(\overline{\theta},\varepsilon)} S_n(\theta)\,-S_n(\overline{\theta})\Bigr)\ =\ 
\liminf_{n\to\infty}\frac{1}{n}\inf_{\theta\in C(\overline{\theta},\varepsilon)} D_n(\theta,\overline{\theta})
\ >0\ \mbox{ a.s.}
\]
Hence \ $\inf_{\theta\in C(\overline{\theta},\varepsilon)} S_n(\theta)\,-S_n(\overline{\theta})\,\asto\infty$.
\eop
 
For deriving asymptotic normality of the adaptive least squares estimators further assumptions are needed.
Firstly, a condition (B5) on the family of functions $f_\theta$, $\theta\in\Theta$,  and the mean response $\mu$ 
is added to conditions (B1)-(B4). Secondly, two additional conditions (L) and (AH) on the error variables 
in (A1)-(A2) are imposed, where `L' stands for `Lindeberg' and `AH' for `asymptotic homogeneity'. 
\begin{itemize}
\item[(B5)] 
$\Theta\subseteq\mathbb{R}^p$ endowed with the usual Euclidean metric, 
${\rm int}(\Theta)\not=\emptyset$, where ${\rm int}(\Theta)$ denotes the interior of $\Theta$ as a subset of $\mathbb{R}^p$,
the function \ $\theta\mapsto\mu(x,\theta)$ \ is twice continuously differentiable on the interior of $\Theta$
for each fixed $x\in{\cal X}$, and  
\[
f_\theta(x)\,=\,\nabla \mu(x,\theta)\quad\mbox{for all $\theta\in{\rm int}(\Theta)$ and all $x\in{\cal X}$,}
\]
where $\nabla \mu(x,\theta)=
\Bigl(\frac{\partial}{\partial\theta_1}\mu(x,\theta),\ldots,\frac{\partial}{\partial\theta_p}\mu(x,\theta)\Bigr)^\trp$
for $\theta=(\theta_1,\ldots,\theta_p)^\trp\in{\rm int}(\Theta)$.
\end{itemize}
For a subset $A\subseteq\Omega$ we denote by $\Ifkt(A)$ the function on $\Omega$ which is constantly equal to $1$ on $A$
and is constantly equal to $0$ on $\Omega\setminus A$. 
\begin{itemize}
\item[(L)]
\hspace*{2ex}$\displaystyle\frac{1}{n}\sum_{i=1}^n {\rm E}\Bigl(e_i^2\Ifkt\bigl(|e_i|>\varepsilon\sqrt{n}\bigr)\,\Big|{\cal F}_{i-1}\Bigr)
\,\asto0$ \ for all $\varepsilon>0$.
\end{itemize}
\begin{itemize}
\item[(AH)]
\hspace*{2ex}$\displaystyle{\rm E}\bigl(e_n^2\big|{\cal F}_{n-1}\bigr)\asto\sigma^2(\overline{\theta})$ \ for some positive real constant
$\sigma^2(\overline{\theta})$.
\end{itemize}
The following two conditions (L') and (L'') are less technical than the Lindeberg condition (L), and each of them implies (L).
\begin{itemize}
\item[(L')]
$\sup_{i\in\mathbb{N}}{\rm E}\bigl(|e_i|^\alpha\big|{\cal F}_{i-1}\bigr)\,<\,\infty$ \ a.s. \ for some real $\alpha>2$.
\item[(L'')] 
The random variables $e_i$, $i\in\mathbb{N}$, are identically distributed, and \ 
$e_i$, ${\cal F}_{i-1}$ \ are independent for each $i\in\mathbb{N}$.
\end{itemize}
In fact, from (L'), observing the trivial inequality \ 
$e_i^2\Ifkt\bigl(|e_i|>\varepsilon\sqrt{n}\bigr)\le |e_i|^\alpha/(\varepsilon\sqrt{n})^{\alpha-2}$,
it follows that
\[
\frac{1}{n}\sum_{i=1}^n {\rm E}\Bigl(e_i^2\Ifkt\bigl(|e_i|>\varepsilon\sqrt{n}\bigr)\,\Big|{\cal F}_{i-1}\Bigr)\,
\le\,\frac{1}{(\varepsilon\sqrt{n})^{\alpha-2}}\,\sup_{i\in\mathbb{N}}{\rm E}\bigl(|e_i|^\alpha\big|{\cal F}_{i-1}\bigr)  
\ \asto0.
\]
From (L'') it follows for all $i\in\mathbb{N}$
\[
{\rm E}\Bigl(e_i^2\Ifkt\bigl(|e_i|>\varepsilon\sqrt{n}\bigr)\,\big|{\cal F}_{i-1}\Bigr)=
{\rm E}\Bigl(e_i^2\Ifkt\bigl(|e_i|>\varepsilon\sqrt{n}\bigr)\Bigr)=
{\rm E}\Bigl(e_1^2\Ifkt\bigl(|e_1|>\varepsilon\sqrt{n}\bigr)\Bigr) \  \mbox{ a.s.}
\]
Hence
\[
\frac{1}{n}\sum_{i=1}^n {\rm E}\Bigl(e_i^2\Ifkt\bigl(|e_i|>\varepsilon\sqrt{n}\bigr)\,\Big|{\cal F}_{i-1}\Bigr)
= {\rm E}\Bigl(e_1^2\Ifkt\bigl(|e_1|>\varepsilon\sqrt{n}\bigr)\Bigr)\ \mbox{a.s.}
\]
and the expectation on the r.h.s.\ converges to zero as $n\to\infty$. Note also that (L'')
implies ${\rm E}\bigl(e_i^2\big|{\cal F}_{i-1}\bigr)={\rm E}\bigl(e_1^2\bigr)=\sigma^2(\overline{\theta})$, say.
Excluding the trivial case $\sigma^2(\overline{\theta})=0$, we see that condition (L'') also implies condition (AH).

\vspace*{1ex}\noindent
{\bf Remark.} Condition (L') was employed by Lai and Wei \cite{Lai-Wei}, Theorem 1 of that paper, and 
by Chen, Hu, and Ying \cite{Chen-Hu-Ying}, condition (C4) on p.~1161 of that paper. 
Condition (L'') meets the assumption of i.i.d.\ error variables of Pronzato \cite{Pronzato} 
for a particular choice of the sequence of sub-sigma-fields ${\cal F}_i$, $i\in\mathbb{N}_0$.
\eop

The $k$-dimensional normal distribution with expectation $0$ and covariance matrix $C$ is denoted by 
${\rm N}(0,C)$, where $C$ is a positive definite $k\times k$ matrix. 
In particular, ${\rm N}(0,I_k)$ is  the $k$-dimensional standard normal distribution, where $I_k$
denotes the $k\times k$ identity matrix.  
For a sequence $W_n$ of $\mathbb{R}^k$-valued random variables, 
convergence in distribution of $W_n$ (as $n\to \infty$) to a $k$-dimensional normal distribution ${\rm N}(0,C)$ 
is abbreviated by $W_n\dto{\rm N}(0,C)$.   
In the following theorem asymptotic normality of the adaptive least squares estimators 
$\widehat{\theta}_n^{\rm\scriptsize (LS)}$ is established.
To some extent our proof is similar to that of Theorem 2 in Pronzato \cite{Pronzato-b}, though the assumptions are different.
Note that, by our Theorem \ref{theo3-1}, the assumption of strong consistency of the adaptive estimators $\widehat{\theta}_n$ employed by the algorithm
is met under (SI) and if $\widehat{\theta}_n=\widehat{\theta}_n^{\rm\scriptsize (LS)}$, $n\ge n_{\rm\scriptsize st}$.   

\begin{theorem}\quad\label{theo3-2}\\
Assume conditions (SI), (B5), (L), and (AH). Moreover, let
the sequence $\widehat{\theta}_n$ of adaptive estimators employed by the algorithm be strongly consistent,
i.e., $\widehat{\theta}_n\asto\overline{\theta}$, and let $\overline{\theta}\in{\rm int}(\Theta)$.  
Then:
\[
\sqrt{n}\,\sigma^{-1}(\thb)\,
M^{1/2}\bigl(\xi_n,\widehat{\theta}_n^{\rm\scriptsize (LS)}\bigr)\,\bigl(\widehat{\theta}_n^{\rm\scriptsize (LS)}-\overline{\theta}\bigr)
\,\dto\,{\rm N}(0,I_p).
\]
Also, denoting by $M_*=M\bigl(\xi^*_{\overline{\theta}},\overline{\theta})$ the information matrix of a locally D-optimal  
design at $\overline{\theta}$, one has
\[
\sqrt{n}\,\bigl(\widehat{\theta}_n^{\rm\scriptsize (LS)}-\overline{\theta}\bigr)
\,\dto\,{\rm N}\bigl(0,\sigma^2(\thb)\,M_*^{-1}\bigr).
\]
\end{theorem}

\vspace*{1ex}\noindent
{\bf Proof.} \ Choose a compact ball $\overline{B}$ centered at $\overline{\theta}$ and such that
$\overline{B}\subseteq{\rm int}(\Theta)$. By Theorem \ref{theo3-1} there is a random variable $N$ 
with values in $\mathbb{N}\cup\{\infty\}$ such that $N<\infty$ a.s. \ and
$\widehat{\theta}_n^{\rm\scriptsize (LS)}\in\overline{B}$ on $\{N\le n\}$ for all integers $n\ge\nst$.
Note that, since $N$ is almost surely finite, $\Ifkt(N\le n)\asto1$ as $n\to\infty$.
Recall our notation introduced earlier: $S_n(\theta)=\sum_{i=1}^n\bigl(Y_i-\mu(X_i,\theta)\bigr)^2$, $n\ge\nst$, $\theta\in\Theta$.
For the gradients of $S_n(\theta)$ w.r.t. $\theta$ one obtains, using (B5),
\begin{equation}
\nabla S_n(\theta)\,=\,-2\sum_{i=1}^n\bigl(Y_i-\mu(X_i,\theta)\bigr)\,\nabla \mu(X_i,\theta),\ \ \theta\in{\rm int}(\Theta).
\label{eq3-2}
\end{equation}
On $\{N\le n\}$ the gradient at $\widehat{\theta}_n^{\rm\scriptsize (LS)}$ is equal to zero, and hence
$\nabla S_n\bigl( \widehat{\theta}_n^{\rm\scriptsize (LS)}\bigr)-\nabla S_n(\overline{\theta})=-\nabla S_n(\overline{\theta})$.
That equation yields, inserting from (\ref{eq3-2})  and $Y_i=\mu(X_i,\overline{\theta})+e_i$ from (A2),
along with some slight manipulations,
\begin{eqnarray}
&& \sum_{i=1}^ne_i\nabla\mu(X_i,\thb)
\,=\,\sum_{i=1}^n\bigl[\mu(X_i,\LSn)-\mu(X_i,\thb)\bigr]\,\nabla\mu(X_i,\LSn)\nonumber\\
&&\phantom{xxxxx}-\,\sum_{i=1}^ne_i\,\bigl[\nabla\mu(X_i,\LSn)-\nabla\mu(X_i,\thb)\bigr]
\qquad \mbox{on $\{N\le n\}$.}\label{eq3-3}
\end{eqnarray}
With $\sigma^2(\thb)$ according to condition (AH) and the locally D-optimal information matrix $M_*$
introduced in the theorem, we firstly show that
\begin{equation}
n^{-1/2}\sigma^{-1}(\thb)\,M_*^{-1/2}\sum_{i=1}^ne_i\nabla\mu(X_i,\thb)\dto{\rm N}(0,I_p).\label{eq3-4}
\end{equation}
To this end, according to the Cram\'{e}r-Wold device, let $v\in\mathbb{R}^p$, $v^\trp v=1$, be given.
Denote \ $Z_i:=\sigma^{-1}(\thb)\,v^\trp M_*^{-1/2}\nabla\mu(X_i,\thb)$ and \ 
$\widetilde{e}_i:=e_iZ_i$, $i\in\mathbb{N}$. Abbreviating 
the random variables on the l.h.s.\ of (\ref{eq3-4}) by $W_n$, one has  $v^\trp W_n\,=\,n^{-1/2}\sum_{i=1}^n\widetilde{e}_i$.
The random variable
$Z_i$ is ${\cal F}_{i-1}$-measurable for all $i\in\mathbb{N}$, and the  
$Z_i$, $i\in\mathbb{N}$, are uniformly bounded: $|Z_i|\le c$ for all $i\in\mathbb{N}$ 
for some positive real constant $c$.
Hence the sequence of partial sums $\sum_{i=1}^n\widetilde{e}_i$, is a martingale w.r.t. ${\cal F}_n$, $n\in\mathbb{N}$,
and we can apply Corollary 3.1 of Hall and Heyde \cite{Hall-Heyde} which states that the following two conditions 
(a) and (b) together imply the distributional convergence \ $n^{-1/2}\sum_{i=1}^n\widetilde{e}_i\dto{\rm N}(0,1)$. \\[.5ex]
(a) $\displaystyle\frac{1}{n}\sum_{i=1}^n{\rm E}\bigl(\widetilde{e}_i^2\big|\,{\cal F}_{i-1}\bigr)\,\asto 1$, \ 
(b) $\displaystyle\frac{1}{n}
\sum_{i=1}^n{\rm E}\Bigl(\widetilde{e}_i^2\Ifkt\bigl(|\widetilde{e}_i|>\varepsilon\sqrt{n}\,\bigr)\big|\,{\cal F}_{i-1}\Bigr)
\,\asto 0$ \ for all $\varepsilon>0$.\\[.5ex]
Condition (b) follows from condition (L) since 
\[
{\rm E}\Bigl(\widetilde{e}_i^2\Ifkt\bigl(|\widetilde{e}_i|>\varepsilon\sqrt{n}\bigr)\,\big|\,{\cal F}_{i-1}\Bigr)
\le c^2{\rm E}\Bigl(e_i^2\Ifkt\bigl(|e_i|>(\varepsilon/c)\sqrt{n}\bigr)\,\big|\,{\cal F}_{i-1}\Bigr).
\]
To verify (a) we write
\begin{eqnarray*}
&&\frac{1}{n}\sum_{i=1}^n{\rm E}\bigl(\widetilde{e}_i^2\big|{\cal F}_{i-1}\bigr)\,=\,
\frac{1}{n}\sum_{i=1}^n{\rm E}\bigl(e_i^2\big|{\cal F}_{i-1}\bigr)\,Z_i^2\\
&&=\,\frac{1}{n}\sum_{i=1}^n\bigl[{\rm E}\bigl(e_i^2\big|{\cal F}_{i-1}\bigr)-\sigma^2(\thb)\bigr]\,Z_i^2
\,+\,\sigma^2(\thb)\frac{1}{n}\sum_{i=1}^nZ_i^2.
\end{eqnarray*}
By (AH) and $|Z_n|\le c$ for all $n\in\mathbb{N}$ one has 
$\bigl[{\rm E}\bigl(e_n^2\big|{\cal F}_{n-1}\bigr)-\sigma^2(\thb)\bigr]\,Z_n^2\asto0$ and hence  
$\frac{1}{n}\sum_{i=1}^n\bigl[{\rm E}\bigl(e_i^2\big|{\cal F}_{i-1}\bigr)-\sigma^2(\thb)\bigr]\,Z_i^2\asto0$.
By the definition of $Z_i$, $i\in\mathbb{N}$, and by (B5),
\begin{eqnarray*}
&&\sigma^2(\thb)\frac{1}{n}\sum_{i=1}^nZ_i^2\,=\,
v^\trp M_*^{-1/2}\Bigl[\frac{1}{n}\sum_{i=1}^n\nabla \mu(X_i,\thb)\,\nabla^\trp\mu(X_i,\thb)\Bigr]M_*^{-1/2}v\\
&&=\,v^\trp M_*^{-1/2}M(\xi_n,\thb)\,M_*^{-1/2}v\,\asto1,
\end{eqnarray*}
where the final convergence is implied by $M(\xi_n,\overline{\theta})\asto M_*$ according to Corollary 3.2 in 
\cite{FF-NG-RS-18}. This proves (a) and hence (\ref{eq3-4}). Next we show that
\begin{eqnarray}
&&n^{-1/2}\sum_{i=1}^n\bigl[\mu(X_i,\LSn)-\mu(X_i,\thb)\bigr]\,\nabla\mu(X_i,\LSn)\nonumber\\
&&\phantom{xxxxx}=\,
\bigl[M(\xi_n,\LSn)\,+\,A_n\bigr]\,\bigl[n^{1/2}\bigl(\LSn-\thb\bigr)\bigr],\label{eq3-5}\\
&&\mbox{with a sequence $A_n$, $n\ge\nst$, of random $p\times p$ matrices such that $A_n\asto0$.}\nonumber
\end{eqnarray}
By the mean value theorem, for each $n$ there are (random) points $\widetilde{\theta}_{i,n}$, $1\le i\le n$,
on the line segment joining $\LSn$ and $\thb$ such that
\[
\mu(X_i,\LSn)-\mu(X_i,\thb)\,=\,\nabla^\trp \mu(X_i,\widetilde{\theta}_{i,n})\,\bigl(\LSn-\thb\bigr).
\]
So we can write, again using (B5),
\begin{eqnarray*}
&&n^{-1/2}\sum_{i=1}^n\bigl[\mu(X_i,\LSn)-\mu(X_i,\thb)\bigr]\,\nabla\mu(X_i,\LSn)\\
&&=\,
n^{-1}\sum_{i=1}^n\nabla\mu(X_i,\LSn)\,\nabla^\trp \mu(X_i,\widetilde{\theta}_{i,n})\,
\bigl[n^{1/2}\bigl(\LSn-\thb\bigr)\bigr]\\
&&=\,\Bigl(M(\xi_n,\LSn)\,+\,\frac{1}{n}\sum_{i=1}^n\nabla\mu(X_i,\LSn)\,
\bigl[\nabla\mu(X_i,\widetilde{\theta}_{i,n})-\nabla\mu(X_i,\LSn)\bigr]^\trp\Bigr)\,\bigl[n^{1/2}(\LSn-\thb)\bigr].
\end{eqnarray*}
For $A_n:=\frac{1}{n}\sum_{i=1}^n\nabla\mu(X_i,\LSn)\,
\bigl[\nabla\mu(X_i,\widetilde{\theta}_{i,n})-\nabla\mu(X_i,\LSn)\bigr]^\trp$ we get, using the Frobenius norm
in the space of  $p\times p$ matrices, i.e., $\Vert A\Vert_{\rm\scriptsize F} = \bigl[{\rm trace}(AA^\trp)\bigr]^{1/2}$,
\begin{eqnarray*}
&&\Vert A_n\Vert_{\rm\scriptsize F}\,\le\,\frac{1}{n}\sum_{i=1}^n\big\Vert \nabla\mu(X_i,\LSn)\,
\bigl[\nabla\mu(X_i,\widetilde{\theta}_{i,n})-\nabla\mu(X_i,\LSn)\bigr]^\trp\big\Vert_{\rm\scriptsize F}\\
&&=\frac{1}{n}\sum_{i=1}^n \big\Vert\nabla\mu(X_i,\LSn)\big\Vert\cdot
\big\Vert\nabla\mu(X_i,\widetilde{\theta}_{i,n})-\nabla\mu(X_i,\LSn)\big\Vert,
\end{eqnarray*}
where we have used that \ $\big\Vert vw^\trp\big\Vert_{\rm\scriptsize F}=\Vert v\Vert\cdot\Vert w\Vert$ for
$v,w\in\mathbb{R}^p$.  
By compactness of ${\cal X}\times \overline{B}$ and uniform continuity 
of $\nabla\mu(x,\theta)$ on ${\cal X}\times \overline{B}$, one has \ 
\[
\overline{c}\,:=\,
\sup_{\theta\in\overline{B},\,x\in{\cal X}}\big\Vert\nabla\mu(x,\theta)\big\Vert
\,<\infty.
\]  
From \ $\max_{1\le i\le n}\Vert\LSn-\widetilde{\theta}_{i,n}\Vert\le\Vert\LSn-\thb\Vert\asto0$
(as $n\to0$) 
and, again, by the uniform continuity of $\nabla\mu(x,\theta)$ on ${\cal X}\times \overline{B}$, one gets 
\[
\Vert A_n\Vert_{\rm\scriptsize F}\le \overline{c}\,
\max_{1\le i\le n}\big\Vert\nabla\mu(X_i,\widetilde{\theta}_{i,n})-\nabla\mu(X_i,\LSn)\big\Vert\,\asto0
\]
which proves (\ref{eq3-5}). Next we show that 
\begin{eqnarray}
&&n^{-1/2}\sum_{i=1}^ne_i\bigl[\nabla\mu(X_i,\LSn)-\nabla\mu(X_i,\thb)\bigr]\,=\,
B_n\,\bigl[n^{-1/2}\bigl(\LSn-\thb\bigr)\bigr],\label{eq3-6}\\
&&\mbox{with a sequence $B_n$ of $p\times p$ random matrices such that $B_n\asto0$.}\nonumber
\end{eqnarray}
Let $v\in\mathbb{R}^p$ be arbitrarily given. We can write, applying the mean value theorem,
\begin{eqnarray}
&&v^\trp\Bigl(n^{-1/2}\sum_{i=1}^ne_i\bigl[\nabla\mu(X_i,\LSn)-\nabla\mu(X_i,\thb)\bigr]\Bigr)\nonumber\\
&&=\,n^{-1/2}\sum_{i=1}^ne_i\bigl[v^\trp\nabla\mu(X_i,\LSn)-v^\trp\nabla\mu(X_i,\thb)\bigr]\nonumber\\
&&=\,n^{-1/2}\sum_{i=1}^ne_iv^\trp \nabla^2\mu\bigl(X_i,\widetilde{\theta}_{i,n}(v)\bigr)\,\bigl(\LSn-\thb\bigr)\nonumber\\
&&=\,\frac{1}{n}\sum_{i=1}^ne_iv^\trp \nabla^2\mu\bigl(X_i,\widetilde{\theta}_{i,n}(v)\bigr)\,\bigl[n^{1/2}\bigl(\LSn-\thb\bigr)\bigr],
\label{eq3-6a}
\end{eqnarray}
where $\nabla^2\mu(x,\theta)$ denotes the Hessian matrix (matrix of second partial derivatives)
of $\mu$ w.r.t. $\theta$ for fixed $x\in{\cal X}$, and $\widetilde{\theta}_{i,n}(v)$, $1\le i\le n$,
are suitable (random) points on the line segment joining $\LSn$ and $\thb$.
Let \ $b_n(v):=\frac{1}{n}\sum_{i=1}^ne_i\nabla^2\mu\bigl(X_i,\widetilde{\theta}_{i,n}(v)\bigr)\,v$ and write
$b_n(v)=b_n^{(1)}(v)+b_n^{(2)}(v)$, where 
\[
b_n^{(1)}(v):=\frac{1}{n}\sum_{i=1}^ne_i\nabla^2\mu(X_i,\thb)\,v\ \mbox{ and }\ 
b_n^{(2)}(v):=\frac{1}{n}\sum_{i=1}^ne_i\bigl[\nabla^2\mu\bigl(X_i,\widetilde{\theta}_{i,n}(v)\bigr)\,v
-\nabla^2\mu(X_i,\thb)\,v\bigr].
\]
Applying Lemma 3.1 (b) in \cite{FF-NG-RS-18} to each component of $b_n^{(1)}(v)$ one gets 
$b_n^{(1)}(v)\asto0$. The uniform continuity of $(x,\theta)\mapsto \nabla^2\mu(x,\theta)\,v$ on ${\cal X}\times\overline{B}$
and $\max_{1\le i\le n}\Vert\widetilde{\theta}_{i,n}(v)-\thb\Vert\le\Vert\LSn-\thb\Vert\asto0$ imply
that $\max_{1\le i\le n}\big\Vert\nabla^2\mu\bigl(X_i,\widetilde{\theta}_{i,n}(v)\bigr)\,v\,-\,\nabla^2\mu(X_i,\thb)\,v\big\Vert\,\asto0$.
By Lemma 3.1 (a) in \cite{FF-NG-RS-18}, $\limsup_{n\to\infty}\frac{1}{n}\sum_{i=1}^n|e_i|\,<\infty$ \,a.s., and hence
\[
\Vert b_n^{(2)}(v)\Vert\le 
\max_{1\le i\le n}\big\Vert \nabla^2\mu\bigl(X_i,\widetilde{\theta}_{i,n}(v)\bigr)\,v\,-\,\nabla^2\mu(X_i,\thb)\,v\big\Vert\,
\frac{1}{n}\sum_{i=1}^n|e_i|\ \asto\,0.
\]
Observing (\ref{eq3-6a}) we have thus obtained that for every $v\in\mathbb{R}^p$ 
\[
v^\trp\Bigl(n^{-1/2}\sum_{i=1}^ne_i\bigl[\nabla\mu(X_i,\LSn)-\nabla\mu(X_i,\thb)\bigr]\Bigr) 
\,=\,b_n^\trp(v)\,\bigl[n^{1/2}\bigl(\LSn-\thb\bigr)\bigr],
\]
where $b_n(v)\asto0$. Specializing to the elementary unit vectors $v^{(\ell)}$, $1\le\ell\le p$, and 
taking the matrix  $B_n$ with rows $b_n^\trp(v^{(\ell)})$, $1\le\ell\le p$, one gets (\ref{eq3-6}).  
So, by (\ref{eq3-3}), (\ref{eq3-4}), (\ref{eq3-5}), and (\ref{eq3-6}) one gets
\[
\sigma^{-1}(\thb)\,M_*^{-1/2}\bigl[M(\xi_n,\LSn)+A_n-B_n\bigr]\,\bigl[\sqrt{n}\bigl(\LSn-\thb\bigr)\bigr]\ 
\dto\,{\rm N}(0,I_p),
\]
where $A_n\asto0$ and $B_n\asto0$. By Theorem \ref{theo3-1} and by Corollary 3.2 of \cite{FF-NG-RS-18},
one has \ $M(\xi_n,\LSn)\asto M_*$. 
Using standard properties of convergence in distribution one gets  
\begin{eqnarray*}
&&\sigma^{-1}(\thb)\,M_*^{1/2}\bigl[\sqrt{n}\bigl(\LSn-\thb\bigr)\bigr]\dto{\rm N}(0,I_p),\quad
\sqrt{n}\bigl(\LSn-\thb\bigr)\dto{\rm N}\bigl(0,\sigma^2(\thb)\,M_*^{-1}\bigr),\\
&&\mbox{and}\quad \sigma^{-1}(\thb)\,M^{1/2}\bigl(\xi_n,\LSn\bigr)\,\bigl[\sqrt{n}\bigl(\LSn-\thb\bigr)\bigr]\dto{\rm N}(0,I_p).
\end{eqnarray*}
\eop


\end{document}